\documentclass[11pt]{amsart}
\usepackage{amsfonts}
\usepackage[all]{xy}
\begin{document}

\newcommand{\half}{{\textstyle{\frac{1}{2}}}}
\newcommand{\aut}{{\rm Aut}}
\newcommand{\et}{{\rm et}}
\newcommand{\A}{{\mathbb A}}
\newcommand{\C}{{\mathbb C}}
\newcommand{\D}{{\mathbb D}}
\newcommand{\F}{{\mathbb F}}
\newcommand{\G}{{\mathbb G}}
\newcommand{\Q}{{\mathbb Q}}
\newcommand{\R}{{\mathbb R}}
\newcommand{\Z}{{\mathbb Z}}
\newcommand{\End}{{\rm End}}
\newcommand{\TAF}{{\rm TAF}}
\newcommand{\Sh}{{\rm Sh}}
\newcommand{\bu}{{\bar{u}}}
\newcommand{\hA}{{\hat{A}}}
\newcommand{\hB}{{\hat{B}}}
\newcommand{\oh}{{\bf o}}
\newcommand{\tF}{{\tilde{F}}}
\newcommand{\Oh}{{\mathbb O}}
\newcommand{\U}{{\rm U}}
\newcommand{\spec}{{\rm Spec}}

\title{Abelian varieties and the Kervaire invariant}
\author{Jack Morava}
\address{Department of Mathematics, Johns Hopkins University, Baltimore,
Maryland 21218}
\email{jack@math.jhu.edu}
\date {1 May 2011 These are notes from a talk at the April 2011 ICMS workshop in Edinburgh. 
It is a pleasure to thank the organizers of that historic meeting.}

\begin{abstract}{This is an entirely expository account of some aspects of the recent solution 
[by Hill, Hopkins, and Ravenel [5]] of the most refractory [7] problem in (higher dimensional, ie 
excluding the Poincar\'e conjecture) differential topology. My focus is the {\bf tools} involved 
in the proof, and their relation to new developments in the theory of topological 
automorphic forms [1,2].}\end{abstract}

\maketitle

\begin{center}{\bf \S 1 Introduction and background} \end{center} \bigskip

\noindent
This began as a presentation in the arithmetic geometry seminar at Johns Hopkins, and is thus
tilted toward number theorists and algebraic geometers. I hope to explain how the Kervaire problem 
involves the 2-adic behavior of certain Shimura varieties, of type $\U(1,3)$. Such moduli stacks 
play a central role in recent work on the local Langlands program [4,6] but the key ideas go back 
thirty years [3]. I owe many thanks to both Andrew Salch and Mark 
Behrens for helpful conversations and correspondence.\bigskip

\noindent
{\bf Some history:}
\bigskip

\noindent
{\bf 1.1} Integers of the form $k = 2^n - 1$ are sometimes called {\bf Mersenne} numbers, but their
study goes back to antiquity: \medskip

\begin{center}
{if $k$ is {\bf prime} then $\half k(k+1) = {2^n \choose 2}$ is {\bf perfect}.}
\end{center}\medskip

\noindent
It will be convenient to allow $n = 0$ and 1 in this definition.\bigskip

\noindent
These numbers appear in topology as relative dimensions of the Hopf fibrations
\[
\xymatrix{
{S^{2d-1} = S(\D^2)} \ar[r]^{h_n} & {\D \cup \infty = S^d} \ar[r] & {\mathbb P}_2(\D)
}
\]
($h(x,y) = x \cdot y^{-1})$ associated to the $d = 2^n$-dimensional division algebras 
$\D = \R,\C,\mathbb{H},\Oh$. \bigskip

\noindent
In {\bf motivic} language, these Hopf maps are represented by classes 
$h_n$ in the $E_2$-term
\[
{\rm Ext}^{1,*}_{\aut(H^*)}(H^*(S^d),H^*(S^{2d-1})) \Rightarrow [S^{2d-1},S^d]^S 
\]
of Adams' [1957] descent spectral sequence for the fiber functor 
$H^*(-,\F_2)$. The group on the right is the $k$th stable homotopy group $\pi^S_k$ (of a point).
\bigskip

\noindent
The classes $h_n$ have analogs in ${\rm Ext}^{1,2^n}$ for all $n$, but when $n > 3$ those analogs 
are killed by Maunder's (1963) differentials 
\[
d_2 \: h_n = h_0 h_{n-1}^2 \;.
\]
Derivatives of $p$th powers vanish identically in characteristic $p$, so this sheds no light
on the behavior of the classes $h_n^2$ in the spectral sequence. \bigskip

\noindent
{\bf 1.2} In one of the foundational papers in differential topology {\bf Kervaire} [1960]
constructed a closed $4m+2$ - dimensional {topological} manifold by `plumbing' together two 
copies of the tangent ball bundle of $S^{2m+1}$, resulting in a space with boundary 
{\bf homeomorphic} to the $(4m+1)$-sphere -- which can thus be capped off.\bigskip

\noindent
He showed that when $m=2$ ({\bf not} a Mersenne number!) the resulting topological manifold could 
possesses {\bf no differentiable structure}; and (using the classification of quadratic forms in 
characteristic two) he identified the obstruction (for general $m$) to smoothability for these manifolds.
\bigskip

\noindent
When $m=0$ (resp. 1) the resulting 2 (resp. 6)-dimensional manifolds are smooth, and it was shown in the 70's
that when $n=4,\; 5$ (ie $k = 15,\;31)$ the resulting manifolds are similarly smoothable; but in 1969 
Browder showed that the $4m+2$-dimensional Kervaire manifold is {\bf not} smoothable {\bf unless} $m$ 
is a Mersenne number, in which case its smoothability is {\bf equivalent} to the survival of the 
corresponding class $h_n^2$ in the Adams spectral sequence\begin{footnote}{Some alternate formulations
of this problem are described in an appendix.}\end{footnote}.\bigskip

\noindent
After about 1980 there was essentially no progress on the question, until the recent work of Hill, 
Hopkins, and Ravenel:\bigskip

\noindent
{\bf 1.3 Theorem:} The Kervaire manifolds are {\bf not} smoothable, {\bf unless} $k$ is
Mersenne with $n \leq 6$ (and perhaps not if $n=6 \; (k = 63)$, which remains open). \bigskip

\noindent
Their arguments show that if $h_n^2$ survives to represent a nontrivial element of 
$\pi^S_{2^{n+1}-2}$, then that class is detected by a certain spectrum $\Omega$ with 
$2^8$-periodic homotopy groups; and that moreover
\[
\Omega_* = 0 \; {\rm if} \; - 4 < * < 0 \; {\rm mod} \; 2^8 \;.
\]
When $* \equiv - 2$ mod $2^8$, ie if $n \geq 7$, this class therefore cannot be nonzero.
\bigskip

\noindent
The proof of this vanishing theorem constitutes the bulk of their paper, which introduces a
host of new techniques (equivariant slice filtrations) in homotopy theory; but this talk
is about questions closer to number theory, which lie at the start of their proof.\bigskip

\noindent
HHR's spectrum $\Omega$ is the (homotopy) fixed point set of a spectrum $\Omega_\Oh$ with 
$\Z_8$-action, and their detection theorem involves the study of the Hurewicz homomorphism
\[
\pi^S_* \to \Omega_* := \Omega_{\Oh^ *}^{h\Z_8}
\]
and its representation
\[
E^{2,*} \to H^2(\Z_8,\Omega^{-*}_\Oh) \to H^2(\Z_8,\Z[\surd i][w^{\pm 1}])
\]
in a descent spectral sequence. This step in the proof follows Ravenel's [1978] work on the $p>3$ 
analog of the Kervaire problem, which involved the study of formal groups of height
$p-1$ at $p$. The recent work is technically much deeper, involving similar questions for formal
groups of height $(p-1)p^2$.\bigskip

\noindent
Here I will try to explain where the spectra $\Omega_\Oh$ and $\Omega$ come from, in the language
of topological modular forms. This was deeply understood by HHR, but the subject is supressed in
their paper, apparently so as not to frighten the horses \dots \bigskip

\begin{center}{\bf \S 2 Moduli of Abelian varieties with prescribed symmetries}\end{center}
\bigskip

\noindent
[The utility of Abelian varieties in this context may have to do with 
their abundant symmetries: they seem to provide natural models for the study of
a large class of orbifold-like singularities \dots] \bigskip

\noindent
{\bf 2.1} The {\bf Shimura stack}
\[
\Sh \; := \; \Sh(p,F,B,V,K\dots) \; \sim \; K \backslash GV(\A)/GV(\Q)
\]
parametrizes Abelian varieties $A$ of dimension $n^2$ with `complex multiplication', ie an
embedding 
\[
\oh_B \to \End(A)
\]
of a maximal order in an algebra $B$ simple over its center $F$, an imaginary quadratic
extension of $\Q$ in which the prime $p = u \bu$ splits: for example, $F = \Q(i), \; p = (1 +
i)(1-i) = 2$, and $B = M_n(F)$.\bigskip

\noindent
Further bells and whistles include \medskip

\noindent
$\; \bullet \;$ a rank one $B$-module $V$, endowed with a $\Q$-valued alternating Hermitian
form [6] of signature $(1,n-1)$, \medskip

\noindent
$\; \bullet \;$ the $\Q$-algebraic group $GV$ of its unitary similitudes, and \medskip

\noindent
$\; \bullet \;$ a sufficiently large open compact subgroup of $GV(\A^{p,\infty})$ (defining a kind
of level structure). \bigskip

\noindent
If $B$ is split (as in the example above) this moduli problem is Morita equivalent to the
classification of $n$-dimensional Abelian varieties with complex multiplication by $F$, in a more
classical sense.\bigskip

\noindent
{\bf 2.2} [1 \S 8.3] {\bf Definition}: $\TAF_{GV}(K)$ is the spectrum defined by global sections 
of a certain sheaf of $E_\infty$ ring-spectra (ie, of `topological automorphic forms') over 
$\Sh(K)\hat{}_{p,\et}$.\bigskip

\noindent
{\bf Theorem} [1 \S 14.5.6]: The completion
\[
\TAF_{GV}(K)_{K(n)} \; \cong \; \prod_{x \in \Sh^{[n]}} EO(x)
\]
(of the Behrens-Lawson spectrum at the $n$th chromatic prime over $p$) is a finite product of
spectra representing certain `higher real $K$-theories' indexed by the 0-dimensional Igusa [4 \S IV]
subscheme $\Sh^{[n]}(K)(\F_p)$ of AV's with `associated' $p$-divisible group whose formal
subgroup has maximal (= $n$) height over $\F_p$.\bigskip

\noindent
This formalism requires considerable unpacking!\bigskip

\noindent
{\bf Examples:} when $n=1$ and $p=2$, $EO(x)$ is the completion of classical real $K$-theory,
interpreted as the homotopy fixed-point spectrum of the action of ${\rm Gal}(\C/\R) = \Z_2$ on
2-adically completed complex $K$-theory. When $n=2$ and $p=2$ or 3, we get topological {\bf modular} 
forms, ie elliptic cohomology, as  in [BL {\tt arXiv0901.1838}].\bigskip

\noindent
{\bf 2.3} The {\bf key idea} [3] is that at $p$, the completion
\[
\oh_{\hB} \cong \oh_{B,u} \times \oh_{B,\bu} \cong M_n(\Z_p)^2 \to \End(\hA) \to \End(A(p))
\]
of complex multiplication {\bf splits} the $p$-divisible group 
\[
A(p) \; \cong \; (\Q/\Z)_p^{2n^2} \; \cong \; A(u) \times A(\bu)
\]
into components of dimension $n$ and $n(n-1)$ respectively: the signature condition forces 
${\rm Lie} \; A \otimes_{\oh_F} \oh_{F,u}$ to be a locally free $n$-dimensional $\oh_{F,u}$-module.
\bigskip

\noindent
The choice of a rank one idempotent $\epsilon \in M_n(\Z_p)$ then defines an isomorphism
\[
A(u) \cong (\epsilon A(u))^n \;,
\]
and $x \in \Sh^{[n]}(\F_p)$ iff the mod $p$ reduction of the formal subgroup of $\epsilon A(u)$ 
has height $n$. In that case
\[
EO(x) = E_n^{h\aut(x)} \;,
\]
where $E_n$ is the Lubin-Tate spectrum which parametrizes liftings of a height $n$ formal
group and $\aut(x)$ is a {\bf finite} group of automorphisms of $x$ (depending on the `level
structure' $K$). When $n=1$, $E_n$ is the $p$-adic completion of complex $K$-theory, and 
$\aut(x)$ is (at $p=2$) just the complex Galois group. More generally,\bigskip

\noindent
{\bf Theorem:} [2 \S 8.1, 9.4]: $EO(x)$, with $\aut(x)$ a {\bf maximal} finite subgroup of
automorphisms of $x$, appears in $\TAF_{GV}$ iff $n = (p-1)p^r$, with $p \in \{2,3,5,7\}$ and 
$r \geq 0$.\bigskip

\noindent
In particular, a point $x \in Sh^{[4]}(\F_2)$ represents the reduction of a formal group
over $R_* = \Z[\surd i][w^{\pm 1}]$ with $\Z_8$ in its automorphism group. The detection
theorem shows that an element of $\pi^S_{2^{n+1}-2}$ representing $h_n^2$ maps to a nontrivial 
class in $H^2(\Z_8;R_*)$. \bigskip

\begin{center}{\bf \S 3 An old question of W. Browder} \end{center} \bigskip

\noindent
From now on I will focus on the case $p=2,\; n=4$, and $F = \Q(i)$; I will also write $\tF = 
\Q(\surd i)$ (which is {\bf wildly} ramified over $\Q$). \bigskip

\noindent
The detecting theory $\Omega$ is then a constituent of $\TAF$, corresponding to a formal group 
\[
x = \epsilon A(u)^0
\]
of height four at $p=2$ with `complex multiplication' by $\Z[\surd i] = \oh_\tF$. More precisely, 
for suitably large $K$, $\TAF$ contains a point with $\aut(x)$ maximal (ie cyclic of order 8). 
The Real structure of $MU$ endows its four-fold smash product with a an action of $\Z_8$ (and 
hence of its group ring $\oh_\tF$).\bigskip

\noindent
The map 
\[
x : A_x(u) \cong (\epsilon A_x(u))^4 \to \spec \; MU^{\wedge 4} H(\Z^4,2) \;,
\]
classifying $x$ is equivariant [HHR \S 11.4] with respect to the natural action of 
\[
M_4(\Z) \otimes \oh_\tF \cong M_4(\oh_\tF)
\]
on the left, and of 
\[
\oh_\tF \otimes M_4(\Z) \cong M_4(\oh_\tF)
\]
on the right.\bigskip

\noindent
From the slice spectral sequence HHR construct [\S 9.11] equivariant classes 
$D \in MU^{\wedge 4}_{19\rho_8}$ and $\Delta_1^{16} \in D^{-1}MU^{\wedge 4}_{2^8}$.
The classifying map $x$ sends $D$ to (a unit times) $w^{19} \in \oh_\tF[w^\pm 1] = R_*$, while 
$\Delta_1^8$ maps to $w^8 \in H^0(\Z_8,R_*)$ and thus factors through the $\Z_8$-spectrum 
$D^{-1}MU^{\wedge 4} \sim \Omega_\Oh$, defining its periodicity operator, yielding a lift
\[
\xymatrix{
\$ \ar[d] \ar[r] & MU^{\wedge 4} \ar[d] \ar[dr] \ar[r]^{\lambda_x} & {\bf taf}_x \ar[d] \\
\Omega = \Omega^{h\Z_8}_\Oh \ar[r] & D^{-1}MU^{\wedge 4} \ar[r]^= & \Omega_\Oh \;.}
\]
[{\bf taf} is in boldface because it's in progress.] One way to interpret the work of HHR is
as first showing that the Kervaire invariant is a topological automorphic form of special 
type; after which they show that in high dimensions it is of {\bf very} special type, ie 
zero. \bigskip

\noindent
Bill Browder asked long ago if one of the root difficulties of the Kervaire invariant problem
might be that $4 = 2 \times 2$ rather than $2 + 2$; this may also be relevant to
Atiyah's ideas about the Freudenthal construction. The work of B \& H suggests that in fact 
\[
4 = (2 - 1) \cdot 2^2 \;.
\]
This is not at all a joke: it's one aspect of a growing consciousness that categorification
can be a powerful tool for revealing structures (eg groups or vector spaces) underlying 
invariants like cardinalities or dimensions. \bigskip

\begin{center}{\bf \S 4 Coda: I'm with Vigleik}\end{center} \bigskip

\noindent
I'd like to close with something I learned from V. Angeltveit (via Andrew Salch): one of the 
first examples of 2-adic exceptionalism is that (if $r \gg 0$) the isomorphism
\[
(\Z/p^r\Z)^\times \; \cong \; \Z/(p-1)p^{r-1}\Z
\]
at odd primes becomes 
\[
(\Z/2^r\Z)^\times \; \cong \; \{\pm 1\} \times \Z/2^{r-2}\Z
\]
at $p=2$.\bigskip

\noindent
VA observes that $\alpha_1(p) \in \pi^S_{2p-3}$, which maps to a nontrivial 
element in the algebraic $K$-theory of $\Z/p^r\Z$, accounts for this anomaly at 
$p=2$ if we believe that
\[
K_1(\Z/p^r\Z) \; \cong \; (\Z/p^r\Z)^\times \;.
\]
At larger primes $\alpha_1$ continues to map to interesting [Soul\'e] classes in $K$-theory; 
when $p = 2^{n-1} + 1$ is small (ie 3 or 5) these are components of the other Hopf maps.  
\bigskip

\noindent
Ravenel's odd-primary analog
\[
d_{2p-1}b_i \equiv h_0b_{i-1}^p \; {\rm mod} \; \dots
\]
of Maunder's differential precludes the existence of odd-primary versions of the Kervaire 
elements; but for all we know, there could be anything from odd perfect numbers to a solution 
to Vandiver's conjecture floating around in the higher reaches of the Adams-Novikov spectral sequence \dots

\newpage

\begin{center}{\bf Appendix: Some digressions} \end{center} \bigskip

\noindent
{\bf 1} Let $M_0, \: M_1$, be closed connected smooth manifolds (of the same dimension)
with embedded top-dimensional disks $D_i \subset M_i$. Let $B(M_i)$ be the 
closed unit tangent ball bundles, and let $B(D_i)$ denote the restrictions of 
those bundles to the chosen disks. \bigskip

\noindent
{\bf Plumbing} the manifolds $M_i$ together along the embedded disks
defines the quotient $P(M_i,D_i)$ of $\coprod B(M_k)$ by the 
identification
\[
B(D_0) \cong D \times D_0 \cong D \times D \to D \times D \cong D \times D_1 \cong B(D_1) \;,
\]
with the arrow being the twist map $(x,y) \mapsto (y,x)$. \bigskip

\noindent
If $m = 0$ as above, both manifolds are circles, $D^1 \times D^1$ is a 
square, and we obtain the topologists' Easter basket: a torus minus a disk.
[Omitting the twist gives a sphere minus two disks, ie the more usual basket.]
\bigskip

\noindent
{\bf 2} The solution to the Kervaire problem has several important (equivalent)
formulations [5]: for example, any stably framed manifold of sufficiently high
(ie $> 126$) dimension is framed cobordant to a homotopy sphere. In purely
homotopy-theoretic terms, the Whitehead square $[\iota,\iota] \in \pi_{2n+1}S^{n+1}$ 
is not divisible by two, except in a finite list of exceptional dimensions. \bigskip

\bibliographystyle{amsplain}

\end{document}